\documentclass[11pt, reqno]{amsart}

\usepackage[OT2, T1]{fontenc}

\usepackage[usenames,dvipsnames]{color}

\usepackage{url}
\usepackage{amsmath}
\usepackage{array}
\usepackage{graphicx}
\usepackage{amssymb}
\usepackage{amsthm}
\definecolor{link}{RGB}{11,0,128}
\usepackage[colorlinks=true,citecolor=NavyBlue,linkcolor=Bittersweet,urlcolor=link]{hyperref}
\usepackage{hyperref}
\usepackage{paralist}
\usepackage{colonequals}		
\usepackage{color}
\usepackage{mathabx}		
\usepackage{amscd}
\usepackage[all,cmtip]{xy}	
\usepackage{sseq}			
\usepackage{verbatim}
\usepackage{parskip}
\usepackage{microtype}
\usepackage[in]{fullpage}
\usepackage{mathrsfs} 			
\usepackage{cleveref}
\usepackage{enumitem}
\usepackage{moreenum}			
\usepackage[alphabetic,lite]{amsrefs} 	
\usepackage{stmaryrd}	
\usepackage[foot]{amsaddr}
\usepackage{subfiles}
\usepackage{ragged2e}

\DeclareSymbolFont{cyrletters}{OT2}{wncyr}{m}{n}
\DeclareMathSymbol{\Sha}{\mathalpha}{cyrletters}{"58}

\newcommand{\bC}{\mathbb{C}}

\newcommand{\bP}{\mathbb{P}}

\newcommand{\bZ}{\mathbb{Z}}

\newcommand{\fm}{\mathfrak{m}}
\newcommand{\fn}{\mathfrak{n}}

\newcommand{\fp}{\mathfrak{p}}

\newcommand{\sE}{\mathscr{E}}

\newcommand{\sH}{\mathscr{H}}

\newcommand{\sL}{\mathscr{L}}

\newcommand{\sO}{\mathscr{O}}

\newcommand{\sV}{\mathscr{V}}

\newcommand{\ra}{\rightarrow}

\newcommand{\xra}{\xrightarrow}

\newcommand{\hra}{\hookrightarrow}

\newcommand{\wt}{\widetilde}
\newcommand{\wh}{\widehat}

\newcommand{\ce}{\colonequals}

\renewcommand{\b}{\textbf}

\newcommand{\isomto}{\overset{\sim}{\longrightarrow}}

\newcommand{\llb}{\llbracket}		
\newcommand{\rrb}{\rrbracket}		
\newcommand{\MCM}{\underline{\mathrm{MCM}}}  

\providecommand{\p}[1]{\left(#1\right)}
\providecommand{\up}[1]{{\upshape(}#1{\upshape)}}
\providecommand{\uref}[1]{{\upshape\ref{#1}}}

\DeclareMathOperator{\Coker}{Coker}		
\DeclareMathOperator{\Spec}{Spec}		
\DeclareMathOperator{\Proj}{Proj}		
\DeclareMathOperator{\Hom}{Hom}			
\DeclareMathOperator{\depth}{depth}		
\DeclareMathOperator{\Ext}{Ext}			
\DeclareMathOperator{\End}{End}		
\DeclareMathOperator{\Pic}{Pic}		
\newcommand{\ba}{\begin{aligned}}
\newcommand{\ea}{\end{aligned}}
\newcommand{\be}{\begin{equation}}
\newcommand{\ee}{\end{equation}}
\newcommand{\pf}{\begin{proof}}
\newcommand{\bpf}{\begin{proof}}
\newcommand{\epf}{\end{proof}}
\newcommand{\bthm}{\begin{thm}}
\newcommand{\ethm}{\end{thm}}

\newcommand{\bprop}{\begin{prop}}
\newcommand{\eprop}{\end{prop}}
\newcommand{\bcor}{\begin{cor}}
\newcommand{\ecor}{\end{cor}}

\newcommand{\brem}{\begin{rem}}
\newcommand{\erem}{\end{rem}}

\newcommand{\brems}{\begin{rems} \hfill \begin{enumerate}[label=\b{\thenumberingbase.},ref=\thenumberingbase]}

\newcommand{\erems}{\end{enumerate} \end{rems}}

\newcommand{\begs}{\begin{egs} \hfill \begin{enumerate}[label=\b{\thenumberingbase.},ref=\thenumberingbase]}

\newcommand{\eegs}{\end{enumerate} \end{egs}}

\newcommand{\blem}{\begin{lemma}}
\newcommand{\elem}{\end{lemma}}

\newcommand{\bconj}{\begin{conj}}
\newcommand{\econj}{\end{conj}}
\newcommand{\bprob}{\begin{Problem}}
\newcommand{\eprob}{\end{Problem}}

\newcommand{\bq}{\begin{Q}}
\newcommand{\eq}{\end{Q}}
\newcommand{\benum}{\begin{enumerate}[label={{\upshape(\alph*)}}]}
\newcommand{\benuma}{\begin{enumerate}[label={{\upshape(\arabic*)}}]}
\newcommand{\benumr}{\begin{enumerate}[label={{\upshape(\roman*)}}]}
\newcommand{\eenum}{\end{enumerate}}
\newcommand{\bc}{}
\newcommand{\bd}{\begin{defn}}
\newcommand{\ed}{\end{defn}}

\newcommand{\beg}{\begin{eg}}
\newcommand{\eeg}{\end{eg}}

\newcommand{\bcl}{\begin{claim}}
\newcommand{\ecl}{\end{claim}}
\newcommand{\lab}{\label}

\newcommand{\x}{\text}

\newcommand{\q}{\quad}
\newcommand{\qq}{\quad\quad}
\newcommand{\qqq}{\quad\quad\quad}

\newcommand{\tst}{\textstyle}

\usepackage{aliascnt}
\newaliascnt{numberingbase}{subsection}

\theoremstyle{plain}
\newtheorem{thm}[numberingbase]{Theorem}
\Crefname{thm}{Theorem}{Theorems}

\Crefname{rethm}{Theorem}{Theorem}
\newtheorem{prop}[numberingbase]{Proposition}
\Crefname{prop}{Proposition}{Propositions} 
\newtheorem{Q}[numberingbase]{Question}
\Crefname{Q}{Question}{Questions}
\newtheorem{Problem}[subsection]{Problem}
\Crefname{Problem}{Problem}{Problems}
\newtheorem{conj}[numberingbase]{Conjecture}
\Crefname{conj}{Conjecture}{Conjectures}
\newtheorem{cor}[numberingbase]{Corollary}
\Crefname{cor}{Corollary}{Corollaries}
\newtheorem{lemma}[numberingbase]{Lemma}

\Crefname{subprop}{Proposition}{Propositions}

\Crefname{subcor}{Corollary}{Corollaries}

\Crefname{sublem}{Lemma}{Lemmas}

\theoremstyle{remark}
\newtheorem{claim}[equation]{Claim}
\Crefname{claim}{Claim}{Claims}

\Crefname{subrem}{Remark}{Remarks}

\theoremstyle{definition}
\newtheorem{defn}[numberingbase]{Definition}
\Crefname{defn}{Definition}{Definitions}

\Crefname{conv}{Convention}{Conventions}
\newtheorem{eg}[numberingbase]{Example}
\Crefname{eg}{Example}{Examples}
\newtheorem{rem}[numberingbase]{Remark}
\Crefname{rem}{Remark}{Remarks}
\newtheorem*{rems}{Remarks}

\theoremstyle{plain}
\newtheorem{thm-tweak}[subsection]{Theorem}
\Crefname{thm-tweak}{Theorem}{Theorems}
\newtheorem{lemma-tweak}[subsection]{Lemma}
\Crefname{lemma-tweak}{Lemma}{Lemmas}
\newtheorem{cor-tweak}[subsection]{Corollary}
\Crefname{cor-tweak}{Corollary}{Corollaries}
\newtheorem{prop-tweak}[subsection]{Proposition}
\Crefname{prop-tweak}{Proposition}{Propositions} 
\newtheorem{conj-tweak}[subsection]{Conjecture}
\Crefname{conj-tweak}{Conjecture}{Conjectures} 

\theoremstyle{definition}
\newtheorem{defn-tweak}[subsection]{Definition}
\Crefname{defn-tweak}{Definition}{Definitions}
\newtheorem{eg-tweak}[subsection]{Example}
\Crefname{eg-tweak}{Example}{Examples}
\newtheorem*{rems-tweak}{Remarks}
\newtheorem{rem-tweak}[subsection]{Remark}
\Crefname{rem-tweak}{Remark}{Remarks}

\newtheoremstyle{subsection-tweak}
   {11pt}
   {3pt}%
   {}
   {}%
   {\bfseries}
   {}%
   {.5em}
   {\thmnumber{\@{#1}{}\@{#2}.}%
    \thmnote{~{\bfseries#3.}}}    
    
\theoremstyle{subsection-tweak}
\newtheorem{pp}[numberingbase]{}
\newcommand{\bpp}{\begin{pp}}
\newcommand{\epp}{\end{pp}}

\theoremstyle{subsection-tweak}
\newtheorem{pp-tweak}[subsection]{}

\numberwithin{equation}{numberingbase}

\makeatletter
\def\@tocline#1#2#3#4#5#6#7{
    \begingroup 
    \@ifempty{#4}{%
    }{%
    }%

    \parindent\z@ \leftskip#3\relax \advance\leftskip\@tempdima\relax
    #5\hskip-\@tempdima
      \ifcase #1
       \or\or \hskip 2em \or \hskip 1em \else \hskip 3em \fi%
      #6\nobreak\relax
    \dotfill\hbox to\@pnumwidth{\@tocpagenum{#7}}\par
    \nobreak
    \endgroup
  }
 \def\l@section{\@tocline{1}{0pt}{1pc}{}{}}

\renewcommand{\tocsection}[3]{%
  \indentlabel{\@ifnotempty{#2}{\makebox[1.3em][l]{%
    \ignorespaces#1 \bfseries{#2}.\hfill}}}\bfseries{#3}
    \vspace{1.5pt}}

\renewcommand{\tocsubsection}[3]{%
  \indentlabel{\@ifnotempty{#2}{\hspace*{-0.5em}\makebox[2.1em][l]{%
    \ignorespaces#1#2.\hfill}}}#3
    \vspace{1.5pt}}

\makeatother

\begin{document}

\title{Grothendieck--Lefschetz for vector bundles}

\author{K\k{e}stutis \v{C}esnavi\v{c}ius}
\address{CNRS, UMR 8628, Laboratoire de Math\'{e}matiques d'Orsay, Universit\'{e} Paris-Saclay, 91405 Orsay, France}
\email{kestutis@math.u-psud.fr}

\date{\today}


\begin{abstract} According to the Grothendieck--Lefschetz theorem from SGA 2, there are no nontrivial line bundles on the punctured spectrum $U_R$ of a local ring $R$ that is a complete intersection of dimension $\ge 4$. Dao conjectured a generalization for vector bundles $\sV$ of arbitrary rank on $U_R$: such a $\sV$ is free if and only if $\depth_R(\End_R(\Gamma(U_R, \sV))) \ge 4$. We use deformation-theoretic techniques to settle Dao's conjecture. We also present examples showing that its assumptions are sharp and draw consequences for splitting of vector bundles on complete intersections in projective space.  \end{abstract}


\maketitle

\hypersetup{
    linktoc=page,     
}

\renewcommand*\contentsname{}
\q\\
\tableofcontents


\section{The conjecture of Dao}

\bpp[The Grothendieck--Lefschetz theorem]
A key result in local commutative algebra, proved by Grothendieck in SGA 2, says that for a Noetherian local ring $R$ that is a complete intersection (in the sense of \S\ref{conv}) of dimension $\ge 4$, every line bundle on the punctured spectrum $U_R$ is trivial, that is, $\Pic(U_R) = 0$ (see \cite{SGA2new}*{XI, 3.13 (ii)}). In contrast, nontrivial vector bundles $\sV$ may exist on $U_R$ even when $R$ is regular. Nevertheless, a conjecture of Hailong Dao \cite{Dao13}*{7.2.2} predicts that 
\be \lab{Dao-conj}
\text{if a vector bundle $\sV$ on $U_R$ satisfies} \q \depth_R(\End_R(\Gamma(U_R, \sV))) \ge 4, \q \text{then $\sV$ is free,}
\ee
in which case the depth in question equals $\dim(R)$. When $\sV$ is a line bundle, $R \isomto \End_R(\Gamma(U_R, \sV))$ (see \Cref{dep-ref}), so the prediction \eqref{Dao-conj} generalizes the Grothendieck--Lefschetz theorem recalled above. The main goal of the present paper is to establish Dao's conjecture in \Cref{Dao-conj-pf} and to deduce the following consequence for vector bundles on complete intersections in projective space.
\epp

\bthm[\Cref{ci-thm}] \lab{ci-thm-ann}
For a field $k$ and a global complete intersection $X \subset \bP^n_k$ of dimension $\ge 3$, a vector bundle $\sE$ on $X$ is a direct sum of powers of $\sO(1)$ if and only if 
\be \lab{CTA-cond}
H^1(X, \sE nd_{\sO_X}(\sV)(i)) = H^2(X, \sE nd_{\sO_X}(\sV)(i)) = 0 \q \x{for every} \q i \in \bZ.
\ee
\ethm

\brem
Results similar to \Cref{ci-thm-ann} were obtained for $X = \bP^n_\bC$ by Luk and Yau  in \cite{LY93}*{Thm.~B}, for $X = \bP^n_k$ by Huneke and Wiegand in \cite{HW97}*{Thm.~5.2}, and for odd-dimensional hypersurfaces of dimension $\ge 3$ by Dao in \cite{Dao13}*{8.3.4}. In these previous results the condition \eqref{CTA-cond} is weaker: one does not assume the vanishing of $H^2$.
\erem

\bpp[The method of proof]
Our argument for \eqref{Dao-conj} is built on the strategy used by Grothendieck for line bundles and rests on the Lefschetz algebraization theorems from SGA 2. More precisely, we begin by using local cohomology to show that the depth assumption implies unobstructed deformations for $\sV$ and then, after replacing $R$ by its completion, use this to lift $\sV$ to the formal completion along a hypersurface of the punctured spectrum of a complete intersection cut out by fewer hypersurfaces. A Lefschetz theorem from SGA 2 allows us to algebraize the lift and, after taking care to retain the depth assumption, we proceed inductively to eventually reduce to regular $R$. To conclude, we use a theorem of Huneke--Wiegand: if $R$ is regular and $\depth_R(\End_R(\Gamma(U_R, \sV))) \ge 3$, then $\sV$ is free. Examples coming from Kn\"{o}rrer periodicity for maximal Cohen--Macaulay modules over local hypersurfaces show that the depth assumption in \eqref{Dao-conj} is optimal, see \S\ref{S4-sharp}.
\epp

\bpp[A previously known case]
When, in addition to $\depth_R(\End_R(\Gamma(U_R, \sV))) \ge 4$, we also have $\depth_R(\Gamma(U_R, \sV)) \ge 3$, the conjecture \eqref{Dao-conj} was established by Dao in \cite{Dao13}*{7.2.3}. In this case, the assumption on $\Gamma(U_R, \sV)$ allows one to transform the depth condition on $\End_R(\Gamma(U_R, \sV))$ into 
\[
\Ext^2_R(\Gamma(U_R, \sV), \Gamma(U_R, \sV)) = 0.
\]
Due to the results of Auslander--Ding--Solberg \cite{ADS93}, the vanishing of this $\Ext^2$ implies that, after replacing $R$ by its completion, the $R$-module $\Gamma(U_R, \sV)$ lifts to a regular ring, and Dao concludes by using the resulting Tor-rigidity of $\Gamma(U_R, \sV)$. In contrast, we bypass any additional hypotheses on $\Gamma(U_R, \sV)$ by deforming over $U_R$ instead of over $R$.
\epp

\bpp[Notation and conventions] \lab{conv}
A Noetherian local ring $(R, \fm)$ is a \emph{complete intersection} if its $\fm$-adic completion is a quotient of a regular local ring by a regular sequence; as is well known, such an $R$ is Cohen--Macaulay. For a local ring $(R, \fm)$, we let $U_R$ denote its \emph{punctured spectrum}:
\[
U_R \ce \Spec(R) \setminus \{ \fm\}.
\]
We use the definition of the (S$_n$) condition given in \cite{EGAIV2}*{5.7.2}: a finite module $M$ over a Noetherian ring $R$ is (S$_n$) if for every prime ideal $\fp \subset R$ one has 
\[
\depth_{R_\fp}(M_\fp) \ge \min(n, \dim(M_\fp)).
\]
We will mostly use this definition when the support of $M$ is $\Spec(R)$, when $\dim(M_\fp) = \dim(R_\fp)$.
\epp

\subsection*{Acknowledgements}
I thank Hailong Dao for very helpful correspondence. I thank Ekaterina Amerik for encouraging me to include \Cref{ci-thm-ann}. 


\section{The Grothendieck--Lefschetz theorem for vector bundles of arbitrary rank}

In order to implement our deformation-theoretic reduction of Dao's conjecture \eqref{Dao-conj} to the case of a regular $R$, we need to show that $\sV$ deforms and that the deformation inherits the depth assumption. The following lemma uses the Lefschetz theorems from \cite{SGA2new} to achieve this.

\blem \lab{pass-to-tilde}
Let $(\wt{R}, \wt{\fm})$ be a complete local ring that is a complete intersection, let $f \in \wt{\fm}$ be a nonzerodivisor, set $R \ce \wt{R}/(f)$, let $\sV$ be a vector bundle on $U_R$, and consider $j \colon U_R \hra \Spec(R)$.
\benum
\item \lab{PTT-a}
If $\dim(R) \ge 3$ and $j_*(\sE nd(\sV))$ is {\upshape(S$_3$)}, then for any lift $\wt{\sV}$ of $\sV$ to a vector bundle on an open neighborhood $\wt{U}$ of the closed subscheme $U_R \subset U_{\wt{R}}$, we have 
\[
\qq \Gamma(\wt{U}, \sE nd(\wt{\sV}))/f\Gamma(\wt{U}, \sE nd(\wt{\sV})) \isomto \Gamma(U_R, \sE nd(\sV)).
\]

\item \lab{PTT-c}
If $\dim(R) \ge 3$ and $j_*(\sE nd(\sV))$ is {\upshape(S$_n$)} with $n \ge 3$, then for any lift $\wt{\sV}$ as in \uref{PTT-a}, the pushforward $\wt{j}_*(\sE nd(\wt{\sV}))$ along $\wt{j} \colon \wt{U} \hra \Spec(\wt{R})$ is also {\upshape(S$_n$)}.

\item \lab{PTT-b}
If $\dim(R) \ge 4$ and $j_*(\sE nd(\sV))$ is {\upshape(S$_4$)}, then a lift $\wt{\sV}$ as in \ref{PTT-a} exists for some $\wt{U}$. 
\eenum
\elem

\bpf
By our assumptions, $R$ is a complete intersection of dimension $\ge 3$, and so are its thickenings
\[
R_n \ce \wt{R}/(f^n) \qq \text{for} \qq n \ge 1.
\]
By the finiteness theorem \cite{SGA2new}*{VIII, 2.3}, 
both $j_*(\sE nd(\sV))$ and $\wt{j}_*(\sE nd(\wt{\sV}))$ are coherent.

\benum
\item
Let $\sV_n \ce \wt{\sV}/f^n$ be the pullback of $\wt{\sV}$ to $U_{R_n}$ (so that $\sV_1 \cong \sV$), and let $\wh{\sV} \cong \varprojlim_n \sV_n$ be the formal $f$-adic completion of $\wt{\sV}$. The formal $f$-adic completion $(\sE nd(\wt{\sV}))\,\wh{\ }$ is then identified with $\varprojlim_n \sE nd(\sV_n)$, and each $\sE nd(\sV_n)$ is a successive extension of copies of $\sE nd(\sV)$. 
Moreover, since the finite $R$-module $\Gamma(U_R, \sE nd(\sV))$ is of depth $\ge 3$, we have
\be \lab{H1-zero}
\qq\ \  H^1(U_R, \sE nd(\sV)) \cong H^2_\fm(R, \Gamma(U_R, \sE nd(\sV))) = 0, \q\text{so also} \q H^1(U_{R_n}, \sE nd(\sV)) = 0
\ee
for every $n > 0$ (see \cite{SGA2new}*{III, 3.3 (iv)}). It follows that 
\[\ba
\tst\qqq \Gamma(U_R, (\sE nd(\wt{\sV}))\,\wh{\ }\,)/f\Gamma(U_R, (\sE nd(\wt{\sV}))\,\wh{\ }\,) &\tst\cong (\varprojlim_n \Gamma(U_R, \sE nd(\sV_n)))/f(\varprojlim_n \Gamma(U_R, \sE nd(\sV_n))) \\ &\tst\cong \varprojlim_n (\Gamma(U_R, \sE nd(\sV_n))/f\Gamma(U_R, \sE nd(\sV_n))) \\ &\cong \Gamma(U_R, \sE nd(\sV)).
\ea\]
To conclude, we use the local Lefschetz theorem \cite{SGA2new}*{X, 2.1 (i)} to obtain
 \[
\qq \Gamma(\wt{U}, \sE nd(\wt{\sV})) \isomto \Gamma(U_R, (\sE nd(\wt{\sV}))\,\wh{\ }\,).
\]

\item
The complement of $\wt{U}$ in $U_{\wt{R}}$ is a union of finitely many closed points of $U_{\wt{R}}$: indeed, the complement of $\wt{U}$ in $\Spec(\wt{R})$ is of the form $\Spec(\wt{R}/I)$ with $(\wt{R}/I)/f(\wt{R}/I)$ Artinian, so $\dim(\wt{R}/I) \le 1$ (see \cite{BouAC}*{VIII.25, Cor.~2~a)}). Thus, since $\wt{R}$ is Cohen--Macaulay and the finite $\wt{R}$-module $\wt{M} \ce \Gamma(\wt{U}, \sE nd(\wt{\sV}))$ is free on $\wt{U}$, we need to show that for every prime $\fp \subset R$ outside $\wt{U}$,
\be \lab{Sn-want}
\qq \depth_{\wt{R}_\fp}(\wt{M}_\fp) \ge \min(n, \dim(\wt{R}_\fp)).
\ee
Scaling by $f$ is injective on $\wt{j}_*(\sE nd(\wt{\sV}))$ because it is so locally over $\wt{U}$, so $f$ is a nonzerodivisor for $\wt{M}$. 
Moreover, by \ref{PTT-a}, the $R$-module $\wt{M}/f\wt{M}$ is identified with $M \ce \Gamma(U_R, \sE nd(\sV))$. Thus, by \cite{EGAIV1}*{0.16.4.10~(i)} and the (S$_n$) assumption on $M$,
\be \lab{cod-mod-x}
\qqq \dim(\wt{R}) - \depth_{\wt{R}}(\wt{M}) = \dim(R) - \depth_{R}(M) \le \dim(R) - \min(n, \dim(R)).
\ee
The inequality \eqref{Sn-want} for $\fp = \wt{\fm}$ follows:
\[
\qq \depth_{\wt{R}}(\wt{M}) \ge \min(n, \dim(R)) + 1 \ge \min(n, \dim(\wt{R})).
\]
Thus, we may assume that $\fp \in \Spec(\wt{R}) \setminus (\wt{U} \cup \wt{\fm})$, so that $\dim(\wt{R}_\fp) = \dim(\wt{R}) - 1 = \dim(R)$. The upper semicontinuity of codepth \cite{EGAIV2}*{6.11.2~(i)} and \eqref{cod-mod-x} then give the desired
\[
\qqq \depth_{\wt{R}_\fp}(\wt{M}_\fp) = \dim(R) - (\dim(\wt{R}_\fp) - \depth_{\wt{R}_\fp}(\wt{M}_\fp)) \ge \min(n, \dim(R)) = \min(n, \dim(\wt{R}_\fp)).
\]
(Alternatively, we could have finished the argument by applying \cite{EGAIV2}*{5.12.2}.)

\item
The coherent $R$-module $\Gamma(U_R, \sE nd(\sV))$ is of depth $\ge 4$, so, analogously to \eqref{H1-zero}, we have
\be \lab{H2-zero}
\qq H^2(U_R, \sE nd(\sV)) \cong H^3_\fm(R, \Gamma(U_R, \sE nd(\sV))) = 0.
\ee
Consequently, since $f \in \wt{R}$ is a nonzerodivisor, there is no obstruction to deforming $\sV$ to $U_{R_2}$ (see, for instance, \cite{Ill05}*{8.5.3 (b)}), to the effect that $\sV$ lifts to a vector bundle $\sV_2$ on $U_{R_2}$. The obstruction to deforming $\sV_2$ to $U_{R_3}$ is again controlled by $H^2(U_R, \sE nd(\sV))$, so $\sV_2$ lifts to a vector bundle $\sV_3$ on $U_{R_3}$. Proceeding in this way, we lift $\sV$ to a vector bundle $\wh{\sV} \ce \varprojlim_n \sV_n$ on the formal $f$-adic completion of $U_{\wt{R}}$. The local Lefschetz theorem \cite{SGA2new}*{X, 2.1 (ii)} then algebraizes $\wh{\sV}$ to a desired $\wt{\sV}$. \qedhere
\eenum
\epf

Geometrically, the depth condition of \eqref{Dao-conj} amounts to the (S$_4$) requirement for $j_*(\sE nd(\sV))$:

\blem \lab{dep-ref}
For a Noetherian local ring $R$ that is of dimension $\ge 2$ and whose completion $\wh{R}$ is {\upshape(S$_2$)}, and for vector bundles $\sV$ and $\sV'$ on $U_R$, 
\[
\text{the $R$-modules} \qq \Gamma(U_R, \sV) \q \text{and} \q \Hom_R(\Gamma(U_R, \sV), \Gamma(U_R, \sV')) \qq \text{are finite and {\upshape(S$_2$)},}
\]
with the associated coherent sheaves $j_*(\sV)$ and $j_*(\sH om(\sV, \sV'))$, where $j \colon U_R \hra \Spec(R)$.
\elem

\bpf
By \cite{EGAIV2}*{5.10.8}, the (S$_2$) assumption implies that $\Spec(\wh{R})$ has no irreducible component of dimension $\le 1$. Thus, since the formation of $j_*(-)$ commutes with the flat base change to $\wh{R}$, the finiteness assertion follows from \cite{SGA2new}*{VIII, 2.3 (ii)$\Leftrightarrow$(iv)} 
(see also \cite{SGA1new}*{VIII,~1.10}). Since $R$ itself is (S$_2$) (see \cite{EGAIV2}*{6.4.1~(i)}), the (S$_2$) assertion and the claim about $j_*(\sV)$ then follow from \cite{EGAIV2}*{5.10.5}. In general, if a finite $R$-module $M$ is of depth $\ge 2$, then so is any $\Hom_R(M', M)$: if $f \in R$ is a nonzerodivisor for $M$, then 
\[
\Hom_R(M', M)/f\Hom_R(M', M) \subset \Hom_R(M', M/fM),
\]
so that any $g \in R$ that is a nonzerodivisor for $M/fM$ is also a nonzerodivisor for 
\[
\Hom_R(M', M)/f\Hom_R(M', M).
\]
In particular, we conclude that $\Hom_R(\Gamma(U_R, \sV), \Gamma(U_R, \sV'))$ is of depth $\ge 2$. Then, by \emph{loc.~cit.},
\[
\Hom_R(\Gamma(U_R, \sV), \Gamma(U_R, \sV')) \isomto \Gamma(R, j_*(\sH om(\sV, \sV'))). \qedhere
\]
\epf

We are ready for the promised extension of the Grothendieck--Lefschetz theorem:

\bthm \lab{Dao-conj-pf}
Let $(R, \fm)$ be a local ring that is a complete intersection of dimension $\ge 4$ and consider the open immersion $j \colon U_R \hra \Spec(R)$. A vector bundle $\sV$ on $U_R$ is free if and only if $j_*(\sE nd(\sV))$ is {\upshape(S$_4$)} \up{that is, if and only if $\depth_R(\End_R(\Gamma(U_R, \sV))) \ge 4$, see Lemma \uref{dep-ref}}.
\ethm

\bpf
By \Cref{dep-ref}, both $j_*(\sV)$ and $j_*(\sE nd(\sV))$ are coherent. If $\sV$ is free, then so is $\sE nd(\sV)$, so that $j_*(\sE nd(\sV))$ is a direct sum of copies of $\sO_{\Spec(R)}$, and hence is (S$_n$) for any $n$ because $R$ is Cohen--Macaulay. For the converse, we  assume that $j_*(\sE nd(\sV))$ is (S$_4$).

To establish the freeness of $\sV$, we will argue that $j_*(\sV)$ is free. Flat base change to $\wh{R}$ commutes with $j_*(-)$, preserves the depth assumption, 
and descends freeness, so we may assume that $R$ is $\fm$-adically complete. Then $R \cong S/(f_1, \dotsc, f_n)$ for a complete regular local ring $(S, \fn)$ and a regular sequence $f_1, \dotsc, f_n \in \fn$. We will argue by induction on $n$, the case $n = 0$ being supplied by \cite{HW97}*{Cor.~2.9}. 

Suppose that $n \ge 1$ and set $\wt{R} \ce S/(f_1, \dotsc, f_{n - 1})$. By \Cref{pass-to-tilde}, the vector bundle $\sV$ lifts to a vector bundle $\wt{\sV}$ defined on some open neighborhood $\wt{U}$ of $U_R$ in $U_{\wt{R}}$ and the pushforward $\wt{j}_*(\sE nd(\wt{\sV}))$ along $\wt{j}\colon \wt{U} \hra \Spec(\wt{R})$ is (S$_4$). We saw in the proof of \Cref{pass-to-tilde}~\ref{PTT-c} that the complement of $\wt{U}$ in $U_{\wt{R}}$ consists of finitely many prime ideals $\fp \subset \wt{R}$ with $\dim(\wt{R}_\fp) = \dim(R)$. The inductive hypothesis applies to the completion of each such $\wt{R}_\fp$ equipped with the restriction of $\wt{\sV}$ to $U_{\wt{R}_\fp}$, to the effect that the restriction of $\wt{j}_*(\wt{\sV})$ to $U_{\wt{R}}$ is a vector bundle. Another application of the inductive assumption, this time to $\wt{R}$ equipped with $(\wt{j}_*(\wt{\sV}))|_{U_{\wt{R}}}$, then proves that $\wt{j}_*(\wt{\sV})$ is free. It follows that $\wt{\sV}$ is free as well, and hence that so is its base change $\sV$.
\epf

\brem
For a further variant of \Cref{Dao-conj-pf}, see \cite{Asg19}*{8.4}.
\erem



\section{The sharpness of the assumptions} \lab{sharpness}

The following examples illustrate the optimality of the assumptions of \Cref{Dao-conj-pf}.

\bpp[The dimension requirement is sharp]
For a field $k$, consider the local ring
\[
R \ce (k[ x, y, z, t ]/(xy - zt))_{(x, y, z, t)}
\]
that is a complete intersection of dimension $3$. We claim that $\Pic(U_R) \cong \bZ$, to the effect that the dimension $\ge 4$ condition of \Cref{Dao-conj-pf} cannot be weakened to $\ge 3$: indeed, for any line bundle $\sL$ on $U_R$, we have $\sO_{U_R} \isomto \sE nd(\sL)$, so $j_*(\sE nd(\sL))$ is (S$_n$) for every $n$, but $\sL$ will need not be $\sO_{U_R}$.

The equation $xy - zt$ cuts out $X \ce \bP^1_k \times \bP^1_k$ sitting in $\bP^3_k$ via its Segre embedding. Since $\Pic(X) \cong \bZ \oplus \bZ$, with the hyperplane class spanning the diagonal copy of $\bZ$ (see \cite{Har77}*{Eg.~II.6.6.2}), we conclude that the Picard group of the punctured spectrum of the local ring of the vertex of the affine cone over $X \subset \bP^3_k$ is $\bZ$ (see \cite{Har77}*{Ex.~II.6.3}). Since this local ring is $R$, we obtain the claimed $\Pic(U_R) \cong \bZ$.
\epp

\bpp[The (S$_4$) requirement is sharp] \lab{S4-sharp}
For every $n \ge 1$ and every algebraically closed field $k$ of characteristic different from $2$, following a suggestion of Hailong Dao, we will construct a \emph{nonfree} finitely generated
 module $M_n$ over the local, $(2n - 1)$-dimensional, complete intersection ring 
\[
R_n \ce k\llb x, y, u_1, v_1, \dotsc, u_{n - 1}, v_{n - 1}\rrb/(xy + u_1v_1 + \dotsb + u_{n - 1}v_{n - 1})
\]
such that $M_n$ is Cohen--Macaulay of depth $2n - 1$ (that is, $M_n$ is ``maximal Cohen--Macaulay'') and the $R_n$-module $\End_{R_n}(M_n)$ is (S$_3$). Since $U_{R_n}$ is regular, the Auslander--Buchsbaum formula will ensure that $M_n$ defines a vector bundle $\sV_n$ on $U_{R_n}$. For $n \ge 2$, the pushforward $(j_n)_*(\sV_n)$ along $j_n \colon U_{R_n} \hra \Spec(R_n)$ will be given by  $M_n$ (see \cite{EGAIV2}*{5.10.5}), so $\sV_n$ will be nonfree but $(j_n)_*(\sE nd(\sV_n))$ will be (S$_3$) (see \Cref{dep-ref}). Thus, for $n \ge 3$ (when $\dim(R_n) \ge 4$), this will show that the (S$_4$) requirement in \Cref{Dao-conj-pf} cannot be weakened to (S$_3$) (even when $j_* (\sV)$ itself is (S$_n$) for every $n$).

For $n = 1$, we set $M_1 \ce k\llb y\rrb$ with $R_1 = k\llb x, y\rrb/(xy)$, so that $M_1$ is a nonfree maximal Cohen--Macaulay $R_1$-module, $\End_{R_1}(M_1)$ is  (S$_3$) (equivalently, (S$_1$)), and $M_1$ admits the free resolution
\[
\dotsc \xra{y} k\llb x, y\rrb/(xy) \xra{x} k\llb x, y\rrb/(xy) \xra{y} k\llb x, y\rrb/(xy) \xra{x} k\llb x, y\rrb/(xy).
\]
This resolution shows that 
\be \lab{Ext-val}
\Ext^{2i - 1}_{R_1}(M_1, M_1) = 0 \q \text{and} \q \Ext^{2i}_{R_1}(M_1, M_1) \cong k \qq \text{for} \qq i \ge 1.
\ee

To construct the remaining $M_n$ from $M_1$, we will use the Kn\"{o}rrer periodicity theorem \cite{Kno87}*{Thm.~3.1}: for every $n \ge 1$, the stable category $\MCM(R_n)$ of maximal Cohen--Macaulay $R_n$-modules\footnote{The objects of $\MCM(R_n)$ are the maximal Cohen--Macaulay $R_n$-modules and the morphisms are given by 
\[
\Hom_{\MCM(R_n)}(M, M') \ce \Hom_{R_n}(M, M')/\{ f \colon M \ra M' \text{ such that $f$ factors through a finite free $R_n$-module}\},
\]
see \cite{Buc87}*{2.1.1 and 4.2.1} for more details.} is equivalent to its counterpart $\MCM(R_{n + 1})$. Explicitly, in terms of matrix factorizations 
\[
\wt{R}_n^a \xra{\varphi} \wt{R}_n^a \xra{\psi} \wt{R}_n^a \qq \text{with} \qq \psi \circ \varphi = \varphi \circ \psi = xy + u_1v_1 + \dotsb + u_{n - 1}v_{n - 1}
\]
 with $\wt{R}_n \ce k\llb x, y, u_1, v_1, \dotsc, u_{n - 1}, v_{n - 1}\rrb$, Kn\"{o}rrer's functor maps the maximal Cohen--Macaulay module $\Coker(\varphi)$ to $\Coker\p{\p{\begin{smallmatrix} u_n & \psi \\ \varphi & -v_n \end{smallmatrix}}}$, where $\p{\begin{smallmatrix} u_n & \psi \\ \varphi & -v_n \end{smallmatrix}}$ is a map in the matrix factorization
\[
\wt{R}_{n+ 1}^a \oplus \wt{R}_{n+ 1}^a \xra{\p{\begin{smallmatrix} u_n & \psi \\ \varphi & -v_n \end{smallmatrix}}} \wt{R}_{n+ 1}^a \oplus \wt{R}_{n+ 1}^a \xra{\p{\begin{smallmatrix} v_n & \psi \\ \varphi & -u_n \end{smallmatrix}}} \wt{R}_{n+ 1}^a \oplus \wt{R}_{n+ 1}^a \qq \text{of} \qq xy + u_1v_1 + \dotsb + u_{n}v_{n}.
\]
By \cite{Buc87}*{4.4.1 (3)}, the category $\MCM(R_n)$ is naturally triangulated, with the translation being given by the syzygy functor $\Coker(\varphi) \mapsto \Coker(\psi)$ (that is, by $(\varphi, \psi) \mapsto (\psi, \varphi)$ on matrix factorizations), which is its own inverse. Thus, the commutativity of the diagram
\[
\xymatrix{
\wt{R}_{n+1}^a \oplus \wt{R}_{n+1}^a \ar[d]_-{\sim}^-{\p{\begin{smallmatrix} 0 & 1  \\ -1 & 0 \end{smallmatrix}}} \ar[rrr]^-{\p{\begin{smallmatrix} u_n & \varphi \\ \psi & -v_n \end{smallmatrix}}} &&& \wt{R}_{n+1}^a \oplus \wt{R}_{n+1}^a \ar[d]_-{\sim}^-{\p{\begin{smallmatrix} 0 & -1  \\ 1 & 0 \end{smallmatrix}}}\ar[rrr]^{\p{\begin{smallmatrix} v_n & \varphi \\ \psi & -u_n \end{smallmatrix}}} &&& \wt{R}_{n+1}^a \oplus \wt{R}_{n+1}^a \ar[d]_-{\sim}^-{\p{\begin{smallmatrix} 0 & 1  \\ -1 & 0 \end{smallmatrix}}} \\
\wt{R}_{n+1}^a \oplus \wt{R}_{n+1}^a \ar[rrr]^-{\p{\begin{smallmatrix} v_n & \psi \\ \varphi & -u_n \end{smallmatrix}}} &&& \wt{R}_{n+1}^a \oplus \wt{R}_{n+1}^a \ar[rrr]^-{\p{\begin{smallmatrix} u_n & \psi \\ \varphi & -v_n \end{smallmatrix}}} &&& \wt{R}_{n+1}^a \oplus \wt{R}_{n+1}^a
}
\]
shows that the Kn\"{o}rrer equivalence commutes with translations.

In summary, the image of $M_1$ under the $(n - 1)$-fold Kn\"{o}rrer equivalence is a maximal Cohen--Macaulay $R_n$-module $M_n$ such that the ``stabilized $\Ext$'s'' defined as in \cite{Buc87}*{6.1.1} by
\[
\underline{\Ext}^i_{R_n}(M, M') \ce \Hom_{\MCM(R_n)}(M, M'[i])
\]
 satisfy
\[
\underline{\Ext}^i_{R_n}(M_n, M_n) \cong \underline{\Ext}^i_{R_1}(M_1, M_1) \qq \text{for every $i$ and $n$.}
\]
Since each $M_n$ is maximal Cohen--Macaulay, \cite{Buc87}*{6.4.1 (i)} ensures that for $i > 0$ the stabilized $\Ext$'s in question agree with their usual nonstable counterparts, so that \eqref{Ext-val} gives 
\[
\Ext^{2i - 1}_{R_n}(M_n, M_n) = 0 \q \text{and} \q \Ext^{2i}_{R_n}(M_n, M_n) \neq 0 \qq \text{for every} \qq n, i \ge 1.
\]
In particular, $\Ext^2_{R_n}(M_n, M_n) \neq 0$, so  each $M_n$ is nonfree. On the other hand, the vanishing of $\Ext^1_{R_n}(M_n, M_n)$ implies that $\End_{R_n}(M_n)$ fits into an $R_n$-module exact sequence
\[
0 \ra \End_{R_n}(M_n) \ra M_n^{\oplus r_1} \ra M_n^{\oplus r_2} \ra M_n^{\oplus r_3} \ra Q \ra 0.
\]
Since $M_n$ is Cohen--Macaulay, it follows that $H^j_\fm(R_n, \End_{R_n}(M_n)) = 0$ for $n \ge 2$ and $j \le 2$ (see \cite{SGA2new}*{III, 3.3}), so that $\End_{R_n}(M_n)$, which is free over $U_{R_n}$, is (S$_3$), as desired. 
\epp

\brem
The dimensions of the rings $R_n$ are odd. Thus, the failure of the freeness of the modules $M_n$ should be contrasted with the following result \cite{Dao13}*{7.2.5}: for a local ring $R$ of \emph{even} dimension $\ge 4$ whose completion is a quotient of an either equicharacteristic or unramified regular local ring by a nonzerodivisor, a vector bundle $\sV$ on $U_R$ is free if and only if $\depth_R(\End_R(\Gamma(U_R, \sV))) \ge 3$.
\erem





\section{Vector bundles on global complete intersections}

We are ready for the promised splitting criterion for vector bundles on global complete intersections.

\bthm \lab{ci-thm}
For a field $k$ and a global complete intersection $X \subset \bP^n_k$ of dimension $d\ge 3$, a vector bundle $\sV$ on $X$ is a direct sum of powers of $\sO(1)$ if and only if 
\be \lab{CT-cond}
H^1(X, \sE nd_{\sO_X}(\sV)(i)) \cong H^2(X, \sE nd_{\sO_X}(\sV)(i)) \cong 0 \q \x{for every} \q i \in \bZ.
\ee
\ethm

\bpf
By the definition of a global complete intersection, $X$ is the $\Proj$ of the graded ring 
\[
R \ce k[x_0, \dotsc, x_n]/(f_1, \dotsc, f_{n - d}) \q \x{for some homogeneous elements} \q f_i \in k[x_0, \dotsc, x_n].
\]
The sequence $f_1, \dotsc, f_{n - d}$ is $k[x_0, \dotsc, x_n]$-regular: indeed, on the local rings of $\bP^n_k$ at the closed points of $X$ this follows from the dimension requirement, and this then implies the same on each prime localization of $k[x_0, \dotsc, x_n]$ due to the gradedness of the kernel of the multiplication by $f_i$ on $k[x_0, \dotsc, x_n]/(f_1, \dotsc, f_{i - 1})$ (the annihilator of this kernel is homogeneous, so if a prime $\fp$ contains it, then so does the prime generated by the homogeneous elements of $\fp$, see also \cite{EGAII}*{2.2.1}). In particular, the local rings of $R$ are complete intersections. In fact, due to the vanishing properties of cohomology of projective space \cite{EGAIII1}*{2.1.13}, the homogeneous coordinate ring of $X$ is $R$:
\be \lab{R-id}
\tst R \cong \bigoplus_{m \ge 0} \Gamma(X, \sO_X(m)), \q \x{compatibly with the gradings.}
\ee
We consider the graded $R$-module $M \ce \bigoplus_{m \in \bZ} \Gamma(X, \sV(m))$, whose associated $\sO_X$-module is $\sV$ (see \cite{EGAII}*{3.4.4}). Letting $\fm  \ce (x_0, \dotsc, x_n) \subset R$ denote the irrelevant ideal, we deduce that  $M$ defines a vector bundle on $\Spec(R) \setminus \{\fm\}$: indeed, on homogeneous localizations of $R$ this follows from $\sV$ being a vector bundle, and, by \cite{EGAII}*{2.2.1}, this then implies the same on the corresponding usual localizations. Moreover, by \cite{EGAIII1}*{1.4.3.2, 2.1.5.2}, the $R$-module $M$ agrees with the pushforward of its restriction to $\Spec(R) \setminus \{ \fm\}$. Thus, \Cref{dep-ref} ensures that the $R$-modules $M$ and $\End_R(M)$ are finite and (S$_2$) and that $\End_R(M)$, like $M$, agrees with the pushforward of its restriction to $\Spec(R) \setminus \{ \fm\}$. The $\sO_X$-module associated to $\End_R(M)$ is $\End_{\sO_X}(\sV)$ (see \cite{EGAII}*{3.2.6}), so \cite{EGAIII1}*{1.4.3.1, 2.1.5.1} supply the identification
\be \lab{coho-id}
\tst H^{j + 1}_\fm(R, \End_R(M)) \cong \bigoplus_{i \in \bZ} H^{j}(X, \sE nd_{\sO_X}(\sV)(i)) \q \x{for every} \q j \ge 1.
\ee
Suppose that $\sV$ is a direct sum of powers of $\sO(1)$, so that $M$ is a direct sum of shifts of the graded $R$-module $\bigoplus_{m\in \bZ} \Gamma(X, \sO_X(m))$. By using the vanishing properties \cite{EGAIII1}*{2.1.13} of the cohomology of the projective space once more, we deduce from \eqref{R-id} that this $R$-module is free of rank $1$. Then $\End_R(M)$ is a free $R$-module and \eqref{coho-id} together with the fact that $\depth(R_\fm) \ge 4$ (obtained from $R_\fm$ being a complete intersection of dimension $\ge 4$) gives \eqref{CT-cond}. 

Conversely, suppose that \eqref{CT-cond} holds. Then \eqref{coho-id} implies that $H^j_\fm(R, \End_R(M)) = 0$ for $j = 2$ or $j = 3$. Since $\End_R(M)$ is (S$_2$), this vanishing also holds for $j \le 1$. It follows that we have $\depth_{R_\fm}((\End_R(M))_\fm) \ge 4$ (see \cite{SGA2new}*{III, 3.3 (iv)}). \Cref{Dao-conj-pf} then implies that $M_\fm$ is $R_\fm$-free, so that $M$ is a finite projective $R$-module. Since $M$ is also graded, the graded Nakayama lemma \cite{BBHR91}*{I, Prop.~1.1 (2)} 
implies that $M$ is generated by a homogeneous lift of any $R/\fm$-basis of $M/\fm M$. Since $X$ is connected, 
the rank of $M$ is constant, and it follows that $M$ is $R$-free as a graded module. In other words, $M$ is isomorphic to a direct sum of shifts of $R$, to the effect that $\sV$ is isomorphic to a direct sum of powers of $\sO(1)$, as desired.
\epf

\end{document}